\documentclass[11pt]{article}
\usepackage[utf8]{inputenc}
\usepackage[english]{babel}
\usepackage{amsmath,amsfonts,amscd,amssymb}
\usepackage[all]{xy}
\usepackage{makeidx}
\usepackage[paper=a4paper,left=30mm,right=20mm,top=25mm,bottom=30mm]{geometry}
    \newenvironment{dedication}
        {\vspace{6ex}\begin{quotation}\begin{center}\begin{em}}
        {\par\end{em}\end{center}\end{quotation}}
\newtheorem{thm}{Theorem}%[section]%,subsection]
\newtheorem{cor}[thm]{Corollary}%[thm]%[section]
\newtheorem{lem}[thm]{Lemma}%[section]
\newtheorem{prop}[thm]{Proposition}%[thm]%[section]
\newtheorem{deft}[thm]{Definition}%[thm]%[section]
\newtheorem{rek}[thm]{Remark}%[section]
\newtheorem{conj}{Conjecture}
\let\noi=\noindent

\def\N{\mathbb{N}} % entiers naturels
\def\Z{\mathbb{Z}} % entiers relatifs
\def\F{\mathbb{F}}

\def\notin{\mbox{$\in$ \hspace{-.8em}/}} % non {\'e}l{\'e}ment de
\makeindex
\title{On the set of bad primes in the study of the Casas--Alvero conjecture}
\author{Daniel Schaub, Univ Angers, CNRS, LAREMA, SFR MATHSTIC\\F-49000 Angers, France
\\email: daniel.schaub@univ-angers.fr \and
Mark Spivakovsky, Univ Paul Sabatier, CNRS, IMT UMR 5219\\ F-31062 Toulouse, France and\\
CNRS, LaSol UMI 2001, UNAM.
\\email: mark.spivakovsky@math.univ-toulouse.fr }

\begin{document}
\maketitle

\begin{dedication}
Dedicated to Professor Osamu Saeki on the occasion of his sixtieth birthday
\end{dedication}

\begin{abstract} The Casas--Alvero conjecture predicts that every univariate polynomial over a field of characteristic zero having a common factor with each of its derivatives $H_i(f)$ is a power of a linear polynomial. One approach to proving the conjecture is to first prove it for polynomials of some small degree $d$, compile a list of bad primes for that degree (namely, those primes $p$ for which the conjecture fails in degree $d$ and characteristic $p$) and then deduce the conjecture for all degrees of the form $dp^\ell$, $\ell\in\N$, where $p$ is a good prime for $d$. In this paper we calculate certain distinguished monomials appearing in the resultant $R(f,H_i(f))$ and obtain a (non-exhaustive) list of bad primes for every degree $d\in\N\setminus\{0\}$.
\end{abstract}

\section{Introduction}

In the year 2001 Eduardo Casas--Alvero published a paper on higher order polar germs of plane curve singularities \cite{C}. His work on polar germs inspired him to make the following conjecture (according to the testimony of Jos\'e Manuel Aroca, E. Casas communicated the problem orally well before 2001).

Let $K$ be a field, $d$ a strictly positive integer and $f = x^d+a_1x^{d-1}+\cdots+a_{d-1}x+a_d$ a monic univariate polynomial of degree $d$ over $K$. Let
$$
H_i(f) = \binom{d}{i}x^{d-i} + \binom{d-1}{i}a_1x^{d-i-1} + \cdots + \binom{i}{i}a_{d-i}
$$
be the $i$-th Hasse derivative of $f$.

\begin{deft}
The polynomial $f$ is said to be a {\bf Casas--Alvero polynomial} if for each $i\in\{1,\ldots,d-1\}$ it has a non-constant common factor with its $i$-th Hasse derivative $H_i(f)$.
\end{deft}

Note that, by definition, a Casas-Alvero polynomial $f$ has a common root with $H_{d-1}(f)$. In particular, if $\text{char}\ K=0$, it has at least one root $b\in K$, regardless of whether or not $K$ is algebraically closed. Making the change of variables $x\rightsquigarrow x-b$, we may assume that $0$ is a root of $f$, in other words, $a_d=0$. In the sequel, we will always make this assumption without mentioning it explicitly.
\begin{conj} {\bf (Casas--Alvero)}
Assume that $\text{char}\ K=0$. If $f \in K[x]$ is a Casas-Alvero polynomial of degree $d$ with $a_d=0$, then $f(x) =x^d$.
\end{conj}

For $i\in\{1,\ldots,d-1\}$, let $R_i = R(f,H_i(f))\in K[a_1,\ldots,a_{d-1}]$ be the resultant of $f$ and $H_i(f)$. The polynomials $f$ and $H_i(f)$ have a common factor if and only if $R_i=0$. Thus $f$ is Casas--Alvero if and only if the point
$(a_1,\ldots,a_{d-1}) \in V(R_1,\dots,R_{d-1}) \subset K^{d-1}$. In those terms the Conjecture can be reformulated as follows:

\begin{conj}
Let $V=V(R_1,\ldots,R_{d-1}) \subset K^{d-1}$. Then $V =\{0\}$. In other words,
\begin{equation} \label{eq:racine}
\sqrt{(R_1,\ldots,R_{d-1})}=(a_1,\ldots,a_{d-1}) \end{equation} 
or, equivalently,
\begin{equation}
\label{eq:rac2} 
a_i^N \in (R_1,\ldots,R_{d-1}) \text{ for all } i \in \{1,\ldots,d-1\} \text{ and some } N \in \N.
\end{equation}
\end{conj}

If $\text{char}\ K=p >0$, the Conjecture is false in general. The simplest counterexample is the polynomial $f(x)= x^{p+1}-x^p$.

\begin{rek}
Let $K \subset K'$ be a field extension. The induced extension
$$
K[a_1,\ldots,a_{d-1}] \subset K'[a_1,\ldots,a_{d-}]
$$
is faithfully flat. Since the polynomials $R_1,\ldots,R_{d-1}$ have coefficients in $\Z$, (\ref{eq:rac2}) holds in\linebreak
$K[a_1,\ldots,a_{d-1}]$ if and only if it holds in $K'[a_1,\ldots,a_{d-1}]$. Hence the truth of the conjecture depends only on the characteristic of $K$ but not on the choice of the field $K$ itself.
\end{rek}

\begin{rek}
Formulae (\ref{eq:racine}) and (\ref{eq:rac2}) can be interpreted in terms of Gröbner bases. Namely, (\ref{eq:racine})  and (\ref{eq:rac2}) are equivalent to saying that for any choice of monomial ordering and Gröbner basis $(f_1,\ldots,f_s)$ of
$(R_1,\ldots,R_{d-1})$, after renumbering the $f_j$, the leading monomial of $f_j$ is a power of $a_j$. 
\end{rek}

We will write CA$_{d,p}$ for the statement ``The Casas-Alvero conjecture holds for polynomials of degree $d$ over fields of characteristic $p$''. 

The following equivalences are known for each $d$ (\cite{DdJ}, \cite{graf}) : 

CA$_{d,0}$ holds $\iff$ CA$_{d,p}$ holds for some prime number $p$ $\iff$ CA$_{d,p}$ holds for all but finitely many primes $p$.

\begin{deft}
 A prime number $p$ is said to be a \textbf{bad prime for} $d$ if CA$_{d,p}$ is false. If $p$ is not a bad prime for $d$, we will say that $p$ is a \textbf{good prime for} $d$.
\end{deft}

\begin{prop} (\cite{graf}, Propositions 2.2 and 2.6)
Take a strictly positive integer $d$, a prime number $p$ and a non-negative integer $\ell$. Assume that CA$_{d,p}$ holds. Then so do CA$_{dp^\ell,p}$ and CA$_{dp^\ell,0}$.
 \end{prop}

\noi This result suggests the following general approach to the problem :

(1) prove the conjecture for a small number $d$;

(2) compile lists of good and bad primes for $d$;

(3) conclude that CA$_{dp^\ell,0}$ holds for all the primes $p$ that are known to be good for $d$.
 
In particular, this shows the importance of knowing which primes are good or bad for a given degree $d$.

The above approach has been carried out up to $d \le 7$ (\cite{CLO}, \cite{frutos}, \cite{C-S}, \cite{DdJ}, \cite{graf}). Some integers cannot be written in the form $dp^\ell$ where $p$ is a good prime for $d$. For example,
$$
12 = 2^2 \cdot3, \ 20= 2^2 \cdot 5,\ 24= 2^3 \cdot 3, \ 28 = 2^2 \cdot7,\ 30= 2 \cdot 3 \cdot 5,\ 36 = 2^2 \cdot 3^2, \ 40 = 2^3 \cdot 5, \ldots
$$

CA$_{12,0}$ has been proved by \cite{CLO} with the aid of a computer, by using a very clever strategy to cut down the computation of resultants and Gröbner basis. Thus the smallest degree $d$ for which CA$_{d,0}$ is not known is $d=20$.
\medskip

\noi In this paper we show that for each $i \in \{1,\ldots,d-1\}$, the monomials $\left(1-\binom{d}{i}\right)^{d-i} a_{d-i}^d$ and
$(-1)^{(d-1)(d-i)}\binom di^{d-1} a_{d-1}^{d-i}a_{d-i}$ appear in the resultant $R_i$ (unless $i=1$ in which case the two monomials are the same and the coefficient is $(1-d)^{d-1}$). Moreover, the monomials $a_{d-i}^d$ are the only pure powers appearing in any of the $R_i$. We then use these facts to compile a (non-exhaustive) list of bad primes for each $d\in\N_{>0}$, namely all the primes $p$ for which there exists $i \in \{1,\ldots,d-1\}$ such that $p\ \left|\ \binom di-1\right.$.
\bigskip

\noi{\bf Acknowledgement.} The fact that the monomial $\left(1-\binom{d}{i}\right)^{d-i} a_{d-i}^d$ appears in $R_i$ and is the only pure power appearing there was first proved by Rosa de Frutos' in her Ph.D. thesis \cite{frutos}, Proposition 2.2.1, page 17.

\section{A list of bad primes}

Unless otherwise specified, from now till the end of this paper we shall regard the $R_i$ as elements of the polynomial ring
$\Z[a_1,\dots,a_{d-1}]$. 

\begin{thm} \label{thm:main} {\bf(\cite{frutos}, Proposition 2.2.1)}
For each $i \in \{1,\ldots,d-1\}$, the monomial\linebreak $(-1)^{d-i}\left(\binom{d}{i}-1\right)^{d-i} a_{d-i}^d$  appears in the resultant $R_i$. Moreover, the monomials $a_{d-i}^d$ are the only pure powers appearing in any of the $R_i$.
\end{thm}

\noi Proof: The polynomial $R_i$ is the resultant of
$$
f=x^d+a_{d-1}x^{d-1}+ \cdots+ a_{d-1}x
$$
and
$$
H_i(f) = \binom{d}{i}x^{d-i} + \binom{d-i-1}{i}a_1x^{d-i-1} + \cdots + \binom{i+1}{i}a_{d-i-1}x + \binom{i}{i}a_{d-i}.
$$
\newpage

\noi{\bf Notation.} For $i,j\in\{1,\ldots,d-1\}$, we denote by $\widetilde a_{ij}$ the element $\binom {d-j} i a_j$.
\medskip

\noi Note that for all $i \in \{1,\ldots,d-1\}$, $\widetilde a_{i,d-i} = a_{d-i}$.
\medskip

The resultant $R_i$ equals the determinant $D(d,i)$ of the following matrix $M(d,i)$:

\[ \left(
 \begin{array}{c|l} \overbrace{
  \begin{array}{ccccc}
1 \ \  &a_1 \ \ &a_2 \ \ &\cdots&a_{d-i-1} \\ 
0 \ \  &1 \ \ &a_1 \ \ &\cdots&a_{d-i-2} \\
\vdots \  \ & \vdots \ \  &&& \\ 
0 \ \  &\cdots \  &0 \ \ &\cdots& 1 \end{array} }^{d-i}
 & \overbrace{
 \begin{array}{lllllll} 
 a_{d-i}&\cdots &a_{d-1} &0 & \cdots& \cdots & 0\\
 a_{d-i-1}&a_{d-i} &\cdots &a_{d-1} &0& \cdots &0 \\
 & & & & &\vdots & \vdots \\
 a_1& a_2& \cdots &a_{d-i}& \cdots &a_{d-1}&0
\end{array} }^{d}
\\ \hline
 \begin{array}{ccccc} 
\binom d i & \widetilde{a}_{i,1} &\widetilde{a}_{i2}&\cdots&\widetilde{a}_{i,d-i-1} \\ 0& \binom d i&\widetilde{a}_{i,1}&\cdots&\widetilde{a}_{i,d-i-2} \\ \vdots& \vdots & & & \vdots \\ \vdots& \vdots &&& \vdots \\ 
\vdots& \vdots &&& \vdots \\ 
\vdots& \vdots &&& \vdots \\  0& 0 &\cdots &\cdots& 0
\end{array}
&
\begin{array}{lllllll} %\hspace{9pt}
 a_{d-i}&\cdots &0 &0 & \cdots& \cdots & 0\\
 \widetilde{a}_{i,d-i-1}& a_{d-i} & 0 & \cdots & & \cdots &0 \\ \widetilde{a}_{i,d-i-2} &\vdots& & &\vdots & \vdots \\
 \vdots & & & & &\vdots & \vdots \\
 \vdots & & & & &\vdots & \vdots \\
  \vdots& & & & &\vdots & \vdots \\
 \cdots & 0& \binom d i &\widetilde{a}_{i,1}& \widetilde{a}_{i,2}&\cdots & a_{d-i}
\end{array}
\end{array}
 \right)
\]
By definition, the determinant $D(d,i)$ of the $(2d-i)\times(2d-i)$ matrix $M(d,i)=(m_{\ell j})$ is
\begin{equation}
\Delta = \sum\limits_{\sigma \in \Sigma_{2d-i}}(-1)^{\epsilon(\sigma)}m_{\sigma(1),1} m_{\sigma(2),2} \cdots m_{\sigma(2d-i),2d-i},\label{eq:determinant}
\end{equation}
where $\Sigma_{2d-i}$ is the group of permutations of $\{1,\ldots,2d-i\}$ and
\begin{eqnarray*}
\epsilon(\sigma)&=&0\quad\text{if }\sigma\text{ is even}\\
&=&1\quad\text{if }\sigma\text{ is odd}. 
\end{eqnarray*}

First of all, note that the last column of $M(d,i)$ has only one non-zero entry that equals $a_{d-i}$. Hence $a_{d-i}\ \left|\ D(d,i)\right.$. In particular no pure power of $a_j$ can appear in $D(d,i)$ for $j \neq d-i$. 

\begin{rek}\label{rem:place}
The entry $a_{d-i}$ appears only in the last $d$ columns of $M(d,i)$: exactly once in each of the last $i$ columns and exactly twice in each of the columns numbered $d-i+1, d-i+2, \ldots,2d-2i$. 
\end{rek}

\noi By inspection of the matrix $M(d,i)$, we see that

(1) a monomial $\omega$ appearing in $D(d,i)$ cannot be divisible by $a_{d-i}^{d+1}$ 

(2) if $a_{d-i}^d\ \not|\ \omega$, then $\omega$ is not a pure power of $a_{d-i}$

(3) if $\left.a_{d-i}^d\ \right|\ \omega$, then in the notation of formula (\ref{eq:determinant}), $\omega = (-1)^{\epsilon(\sigma)} m_{\sigma(1),1}\cdots m_{\sigma(2d-i),2d-i}$ with $\sigma(j)=j$ for $j \in \{2d-2i+1,\ldots, 2d-1\}$ and $\sigma(j) \in \{j,j-d+i\}$ for
$j \in \{d-i+1,\ldots,2d-2i\}$.

The term in (\ref{eq:determinant}) corresponding to $\sigma=\text{Id}$ is the product of the elements on the main diagonal of
$M(d,i)$; this product is equal to $a_{d-i}^d$. There are other choices of $\sigma\in\Sigma_d$ for which the corresponding summand in (\ref{eq:determinant}) is of the form $ca_{d-i}^d$, where $c\in\Z$. Let us group all of these terms together and calculate the total coefficient of $a_{d-i}^d$ in (\ref{eq:determinant})

Now, the coefficient of $a_{d-i}^d$ in $D(d,i)$ is given by the coefficient of $a_{d-i}^{d-i}$ in the expansion of the determinant
$\tilde{D}(d,i)$ of the $(2d-2i)\times (2d-2i)$-matrix $\tilde M(d,i)$ obtained by deleting the last $i$ lines and the last $i$ columns of
$M(d,i)$. 

The matrix $\tilde{M}(d,i)$ is 

\[ \left( \begin{array}{ccccccccc}
 1 & a_1 & a_2 & \ldots & a_{d-i-1} &a_{d-i} & \ldots & a_{d-1} & 0 \\
 0 &1 & a_1 & \ldots & a_{d-i-2}&a_{d-i-1}&a_{d-i} &\ldots &a_{d-1} \\
 \vdots & \vdots &  & &\vdots & & &\vdots&\vdots \\
0 &\ldots& 0 & \ldots &1 & a_1 & a_2 & \ldots & a_{d-i}\\
\binom{d}{i}&\widetilde{a}_{i,1}& \widetilde{a}_{i,2}&\ldots&\widetilde{a}_{i,d-i-1}&  a_{d-i}&0&\ldots&0 \\
0& \binom{d}{i}&\widetilde{a}_{i,1}&\ldots& \widetilde{a}_{i,d-i-2}&\widetilde{a}_{i,d-i-1}& a_{d-i}&0&\ldots 
\\
\vdots& \vdots & & \vdots& \vdots&  & & & \vdots \\
0 & \ldots &0 & \binom{d}{i} & \ldots &\ldots & \widetilde{a}_{i,d-i-1} &a_{d-i}& 0
 \\0 & 0 & \ldots & \ldots& \binom{d}{i} &\ldots &\dots & \widetilde{a}_{i,d-i+1} &a_{d-i}
\end{array} \right)
\]

In this determinant, the first $d-i$ columns do not contain any $a_{d-i}$ and in the last $d-i$ columns, each $a_{d-i}$ appears two times, once in the first $d-i$ rows, once in the last $d-i$ rows. 

In each of the last $d-i$ columns we have to choose one of the two $a_{d-i}$ and delete the rest of the line and the rest of the column to which it belongs. Fix one such choice. The corresponding monomial $(-1)^{\epsilon(\sigma)}m_{\sigma(1),1} m_{\sigma(2),2}\cdots m_{\sigma(2d-2i),2d-2i}$, $\sigma \in \Sigma_{2d-2i}$, satisfies $\sigma(j) \in \{j,j-d+i\}$ for all $j \in \{d-i+1,\ldots,2d-2i\}$. 

Let
\begin{equation} \label{eq:choiceJ}
J=\left\{\left.\sigma(j)\ \right|\ j \in \{d-i+1,\ldots,2d-2i\}\right\}
\end{equation}
and $J^c=\{1,\ldots,2d-2i\} \setminus J$. Write $J=\{j_1,\dots,j_{d-i}\}\subset\{1,\dots,2(d-i)\}$. 
\begin{equation}
\text{For all }\ q,\ell\in\{1,\dots d-i\}\text{ we have }\ j_q-j_\ell \ne d-i.\label{eq:choiceofJ}
\end{equation}
The set $J^c$ has the same property. Note that, conversely, every set
$$
J = \{j_1,\ldots,j_{d-i}\}\subset\{1,\dots,2(d-i)\}
$$
satisfying (\ref{eq:choiceofJ}) has the form (\ref{eq:choiceJ}) for a suitable $\sigma \in \Sigma_{2d-2i}$. 

The coefficient of the term $a_{d-i}^{d-i}$ in the expansion of $\tilde D(d,i)$ corresponding to a given choice of $J$ is the determinant of the matrix $N(d,i,J)$ obtained from the first $(d-i)$ columns of $\tilde
M(d,i)$ by deleting the rows numbered $j_1,\dots j_{d-i-1},j_{d-i}$.

Let
\begin{equation}\label{eq:k}
k=\#\left(J^c\cap\{d-i+1,d-i+2,\dots2(d-i)\}\right).
\end{equation}
There exists a permutation of the rows of $N(d,i,J)$ such that the resulting matrix is an upper triangular matrix with only $1$ and
$\binom di$ on the main diagonal, where $1$ appears $(d-i-k)$ times and $\binom di$ appears $k$ times. Thus the permutation
$\sigma\in\Sigma_{2d-2i}$ is uniquely determined by $J\cap\{d-i+1,\ldots,2d-2i\}$. 

We have
\begin{equation}
\det N(d,i,J)=\pm\binom{d}{i}^k.\label{eq:detN}
\end{equation}
Now, $\sigma$ is the composition of $k$ transpositions $(j,j-d+i)$ for $j \in J^c \cap \{d-i+1,\ldots,2d-2i\}$. Thus
\begin{equation}
\epsilon(\sigma)\equiv k\mod2.\label{eq:signN=k}
\end{equation}
\noi\textbf{Example}
To illustrate the process, let us take
$$
J=\{d-i+1,d-i+2,\dots,2(d-i)-3,d-i-2,2(d-i)-1,2(d-i)\},
$$
which means that we chose all the occurrences of $a_{d-i}$ lying on the main diagonal in the last $(d-i)$ rows of $\tilde M(d,i)$ except in the column number $2(d-i)-2$ in which case we chose the occurrence of $a_{d-i}$ at the place $(d-i-2,2d-2i-2)$. We have $k=1$.

The resulting matrix $N(d,i,J)$ looks like
\[
\left(
\begin{array}{ccccccc}
1&a_1&a_2& \ldots &\ldots &\ldots & a_{d-i-1} \\
0&1&a_1 &\ldots& \ldots & \ldots &a_{d-i-2} \\
\vdots&\vdots & \vdots &  &  & &\vdots\\
0 & \ldots & 0 &1&a_1&a_2&a_3\\
0 & \ldots & \ldots & 0&0&1&a_1\\
0 & \ldots & \ldots& \ldots & 0&0&1\\
0 & \ldots & \ldots & 0& \binom{d}{i}&\widetilde{a}_{i,1} & \widetilde{a}_{i,2}
\end{array}
\right)
\]
To obtain an upper triangular matrix, we have to apply a cyclic permutation to the rows $d-i$, $d-i-1$ and $d-i-2$ and we obtain that the desired determinant is $\pm\binom di$. 
\medskip

Coming back to the proof of the Theorem, for each $k\in\{0,\ldots,d-i\}$, there are $\binom {d-i} k$ choices of $J$ satisfying (\ref{eq:k}). Combining this with (\ref{eq:detN}) and (\ref{eq:signN=k}) and summing over all $k\in\{0,\ldots,d-i\}$, we get that the coefficient of $a_{d-i}^d$ in $R_i$ is
\begin{eqnarray}
\sum_{k=0}^{d-i}(-1)^k\binom{d-i}{k}\binom{d}{i}^k= (-1)^{d-i}\left(\binom{d}{i}-1\right)^{d-i}
\end{eqnarray}
$\Box$

\begin{cor}
Take a prime number $p$ such that there exists $i \in \{1,\ldots,d-1\}$ for which $p\ \left|\ \binom d i -1\right.$. Then CA$_{d,p}$ is false.
\end{cor}

\noi Proof. Assume that char$(K) = p$. By Theorem \ref{thm:main}, no pure power of any of the $a_i$ appears in any of the $R_j\mod p$. Hence the point of $K^{d-1}$ whose $i$-th coordinate is 1 and all of whose other coordinates are zero belongs to
$V(R_1,\ldots,R_{d-1})$. $\Box$\medskip

Using similar arguments, we obtain the following Theorem.

\begin{thm}
For $i\in\{2,\ldots,d-1\}$, the monomial $(-1)^{(d-1)(d-i)}\binom di^{d-1} a_{d-1}^{d-i}a_{d-i}$ appears in the resultant $R_i$.
The term $(-1)^{(d-1)(d-i)}\binom di^{d-1}a_{d-1}^{d-i}a_{d-i}$ is the unique monomial in (\ref{eq:determinant}) of degree $d-i+1$; all the other monomials appearing in (\ref{eq:determinant}) have degree strictly greater than $d-i+1$.
\end{thm}

\noi Proof. By inspection of the matrix $M(d,i)$, we see that the monomial $(-1)^{(d-1)(d-i)}\binom di^{d-1}a_{d-1}^{d-i}a_{d-i}$  appears in the resultant $R_i$: it is the monomial with 
\begin{eqnarray*}
\sigma(j) &=& d-i+j\quad\,\text{ for } \{1,\ldots,d-i\}\\
&=& j-(d-i)\ \text{ for }  j \in \{d-i+1,\ldots,2d-i-1\}\\
&=&2d-i\ \qquad\text{ for } j=2d-i.
\end{eqnarray*}
Moreover, it is the unique monomial $\omega$ of $R_i$ such that $\left.a_{d-1}^{d-i}\ \right|\ \omega$.

Let us prove the second statement of the Theorem. Let $M^\bullet(d,i)$ be the matrix obtained by deleting the last row and the last column of $M(d,i)$. Let $D^\bullet(d,i) = \det M^\bullet(d,i)$. We need to show that all the monomials appearing in $D^\bullet(d,i)$ have order at least $d-i$ and $a_{d-1}^{d-i}$ is the only one of order exactly $d-i$. 

\begin{rek}\label{rek:prop}
For $\ell,j \in \{1,\ldots,2d-i-1\}$, we have $m_{\ell j}\in\N\setminus \{0\}$ if and only if one of the following conditions holds:

(1) $j \in \{d-i+1,\ldots,d-1\}$ and $\ell = j+d-i$

(2) $j \in \{1,\ldots,d-i\}$ and $\ell \in \{j,j+d-i\}$.
\end{rek}
By Remark \ref{rek:prop}, the last $d-i$ columns of $M^\bullet(d,i)$ do not involve any non-zero constant entries. Hence every monomial $\omega = (-1)^{\epsilon(\sigma)}m_{\sigma(1),1} \cdots m_{\sigma(2d-i-1),2d-i-1}$ appearing in $D^\bullet(d,i)$ has degree at least $d-i$. Moreover, assume that $\deg\ \omega = d-i$. Moreover, for $j \in \{1,\ldots,d-1\}$ one of the conditions (1) or (2) of Remark \ref{rek:prop} holds with $\ell = \sigma(j)$. Let
$$
j(\omega) = \min \{j \in \{1,\ldots,d-1\}\ |\ \sigma(j) = j+d-i\}.
$$
\begin{lem}
We have $j(\omega)= 1$.
\end{lem}

\noi Proof of Lemma. Assume that $j(\omega)>1$, aiming for contradiction. By Remark \ref{rek:prop}, 
\begin{equation} \label{eq:prop}
\text{ if } j \in \{d-i+1,\ldots,d-1\},\ \text {then }  \sigma(j) = j+d-i.
\end{equation} Hence $j(\omega) \le d-i+1$. 

Take a $j \in \{ j(\omega)+d-2,\ldots, 2d-i-1 \}$. Then $j>d-1$. From (\ref{eq:prop}) we obtain $\sigma(j)\notin\{2d-i-1,\dots,2d-i-1\}$. By inspection of the matrix $M^\bullet(d,i)$, it follows that
$$
\sigma(j)\in\{ j(\omega)-1,\ldots, d-i\}.
$$
By descending induction on $j$, we obtain 
\begin{equation} \label{eq:prop2}
\sigma(j)=j-d+1 \text{ whenever }  j \in \{ j(\omega)+d-2,\ldots,2d-i-1 \}.
\end{equation}
 
By equation (\ref{eq:prop2}), we have $\sigma(j(\omega)+d-2)= j(\omega)-1$, so $\sigma(j(\omega)-1) \neq j(\omega)-1$. By Remark \ref{rek:prop}, $\sigma(j(\omega)-1) = j(\omega)-1+d-i$. contradicting the definition of $j(\omega)$. The Lemma is proved. $\Box$
\medskip

The Theorem follows from the Lemma by inspection of the matrix $M^\bullet(d,i)$. $\Box$
%The total coefficient of $a_{d-1}^{d-i}a_{d-i}$ in (\ref{eq:determinant}) after grouping the terms is $\binom %di^{d-1}$.   
\bigskip

\end{document}